\documentclass[10pt,reqno]{amsart} 

\usepackage{amssymb,latexsym}
\usepackage{cite} 

\usepackage[height=190mm,width=130mm]{geometry} 
\usepackage{amsfonts}
\usepackage{amsmath}
\usepackage{amscd}
\usepackage{amssymb}
\usepackage{amsthm}
\usepackage{bbm}
\usepackage{bm}
\usepackage{enumerate}
\usepackage{enumitem}
\usepackage{hyperref}
\usepackage{latexsym}
\usepackage{mathabx}
\usepackage{mathdots}
\usepackage{mathtools}

\usepackage{amsmath}
\usepackage{amssymb}
\usepackage{amsthm}
\usepackage{bbm}
\usepackage{enumerate}
\usepackage{latexsym}
\usepackage{mathdots}
\usepackage{geometry}

\usepackage{hyperref}   

\newcommand{\bbar}{\overline{b} }

\newcommand{\C}{\mathbb{C}}

\newcommand{\del}{\partial}

\newcommand{\E}{\mathcal{E}}

\newcommand{\F}{\mathcal{F}}
\newcommand{\I}{\mathcal{I}}

\newcommand{\Ip}{\I_{+}}

\newcommand{\wt}{\textup{wt}}
\newcommand{\Z}{\mathbb{Z}}

\newcommand{\HH}{\mathbb{H}}

\usepackage{amsmath}
\usepackage{amssymb}
\usepackage{amsthm}
\usepackage{bbm}
\usepackage{enumerate}
\usepackage{latexsym}
\usepackage{mathdots}
\usepackage{geometry}

\usepackage{hyperref}   

\theoremstyle{plain}
\newtheorem{theorem}{Theorem}
\newtheorem{conjecture}{Conjecture}

\theoremstyle{definition}

\theoremstyle{remark}

\numberwithin{equation}{section} 

\newcommand{\rmap}{\longrightarrow}

\newcommand{\g}{\ensuremath{\Gamma}}

\newcommand{\ps}{{\raise 1pt\hbox{\tiny (}}}

\newcommand{\pss}{{\raise 1pt\hbox{\tiny [}}}
\newcommand{\pdd}{{\raise 1pt\hbox{\tiny ]}}}
\newcommand{\pd}{{\raise 1pt\hbox{\tiny )}}}

\newcommand{\bs}{{\raise 1pt\hbox{\tiny [}}}
\newcommand{\bd}{{\raise 1pt\hbox{\tiny ]}}}

\def\cross{\mathinner{\mathrel{\raise0.8pt\hbox{$\scriptstyle>$}}
                 \joinrel\mathrel\triangleleft}}

\usepackage{citehack}
\usepackage{amsmath,amsthm}
\usepackage{amsfonts}
\usepackage{amssymb}
\usepackage{longtable}
\usepackage[matrix,arrow,curve]{xy}

\usepackage{lipsum}

\def\W{\mathcal{W}}
\def\K{\mathcal{K}}

\usepackage{stackrel}

\newcommand{\be}{\begin{equation}}
\newcommand{\ee}{\end{equation}}

\newcommand{\nn}{\nonumber \\}

\newcommand{\Y}{\mathcal{Y}}

\newcommand{\nc}{\newcommand}
\nc{\cali}{\mathcal}
\nc{\on}{\operatorname}
\nc{\Wick}{{\mb :}}

\nc{\ddz}{\frac{\partial}{\partial z}}
\nc{\ch}{\mbox{ch}}
\nc{\Oo}{{\cali O}}
\nc{\cond}{|\,}
\nc{\bib}{\bibitem}
\nc{\pone}{\Pro^1}
\nc{\pa}{\partial}
\nc{\arr}{\rightarrow}
\nc{\larr}{\longrightarrow}
\nc{\ket}{\rangle}
\nc{\bra}{\langle}
\nc{\gam}{\bar{\gamma}}
\nc{\q}{\widetilde{Q}}
\nc{\ep}{\epsilon}
\nc{\su}{\widehat{{\mf s}{\mf l}}_2}
\nc{\sw}{{\mf s}{\mf l}}
\nc{\h}{{\mf h}}
\nc{\n}{{\mf n}}
\nc{\ab}{\mf{a}}
\nc{\is}{{\mb i}}
\nc{\js}{{\mb j}}
\nc{\bi}{\bibitem}
\nc{\He}{{\cali H}}
\nc{\inv}{^{-1}}
\nc{\ol}{\overline}
\nc{\wh}{\widehat}
\nc{\dst}{\displaystyle}

\nc{\delt}{\partial_t}
\nc{\ddt}{\frac{\partial}{\partial t}}
\nc{\delx}{\partial_x}
\nc{\mb}{\mathbf}
\nc{\mf}{\mathfrak}

\nc{\mbb}{\mathbb}
\nc{\Ctt}{\C((t))}
\nc{\Ct}{\C[t,t\inv]}

\nc{\un}{\underline}
\nc{\mc}{\mathcal}
\nc{\BB}{{\mc B}}
\nc{\bb}{{\mf b}}
\nc{\kk}{{\mf k}}
\nc{\frob}{\times}
\nc{\sm}{\setminus}
\nc{\Pp}{{\mathbb P}^1}
\nc{\Aa}{{\mc A}}

\nc{\AutO}{\on{Aut}\Oo}
\nc{\AUTO}{\un{\on{Aut}}\Oo}
\nc{\AUTK}{\un{\on{Aut}}\K}
\nc{\Heout}{\He_{\out}}
\nc{\Hetil}{{\widetilde\He}}
\nc{\wb}{\overline}

\nc{\Res}{\on{Res}}
\nc{\pitil}{\Pi}
\nc{\Ctil}{\wt{C}}
\nc{\auto}{\on{Aut} \Oo}
\nc{\phitil}{\wt{\phi}}
\nc{\gz}{\g_{\vec z}}
\nc{\tensorM}{\bigotimes_{i=1}^N{\mathbb M}_i}
\nc{\tensorW}{\bigotimes_{i=1}^N W_{\nu_i,k}}
\nc{\out}{\on{out}}

\nc{\m}{{\mathfrak m}}

\nc{\gx}{\g^0_{\vec x}}

\nc{\hx}{\He^0_{\vec x}}
\nc{\tensorpi}{\pi_{\nu_1,\ldots,\nu_N}^g}
\nc{\Phizw}{\Phi_{\vec w}({\vec z})}
\nc{\Pro}{{\mathbb P}}

\nc{\De}{\Delta}

\nc{\us}{\underset}

\nc{\Ll}{\mc L}
\nc{\dR}{\on{dR}}

\nc{\T}{{\mc T}}

\nc{\Xn}{\overset{\circ}X{}^n} \nc{\Dn}{\overset{\circ}D{}^n}
\nc{\Dxn}{\overset{\circ}D{}^n_x} \nc{\varphitil}{\wt{\varphi}}

\nc{\lf}{{\mf l}}
\nc{\GL}{{}^L G}
\nc{\Vir}{\on{Vir}}

\begin{document}
\title[Relativity of reductive chain complexes of non-abelian simplexes] 
{Relativity of reductive chain complexes of non-abelian simplexes}   
\author{A. Zuevsky}
\address{Institute of Mathematics \\ Czech Academy of Sciences\\ Zitna 25, Prague \\ Czech Republic}

\email{zuevsky@yahoo.com}






\begin{abstract}
Chain total double complexes with reductive differentials 
for non-abelian simplexes 
with associated spaces are considered. 
It is conjectured that corresponding 
relative cohomology 
is equivalent to the coset space of
 vanishing over non-vanishing functionals 
related to differentials of complexes. 
The conjecture is supported by the theorem 
for the case of spaces of correlation functions   
and generalized connections on vertex operator algebra bundles. 
\end{abstract}
\keywords{Cohomology of non-abelian simplexes; correlation functions}
\maketitle
\section{Conflict of Interest Statement}
The author states that: 

1.) The paper does not contain any potential conflicts of interests. 
\section{Data availability statement}
The author confirms  that: 

\medskip 
1.) The paper does not use any datasets. No dataset were generated during and/or analysed 
during the current study. 

2.) The paper includes all data generated or analysed during this study. 

3.) Data sharing not applicable to this article as no datasets were generated or analysed during the current study.

4.) The data of the paper can be shared openly.  
\section{Introduction} 
\label{introduction} 
It is natural to consider non-abelian simplexes with associated spaces \cite{A, BG, PT}
and corresponding cohomology. 
In \cite{BG} spectral sequences for ordinary simplexes with associated functional 
spaces were studied. 
In \cite{BS} it was shown that 
 the Gelfand-Fuks cohomology 
of vector fields on a smooth compact manifold $M$  
 is isomorphic to the singular cohomology 
of the space of continuous cross sections of a certain  
fiber bundle over $M$. 
Passing to a non-abelian simplex setup, one would be interested in 
construction of explicit examples of chain complexes, spectral sequences,  
and relations to geometrical structures of associated manifolds. 
In \cite{BS, Kaw, PT, Fei, Fuks, Wag} cohomology of non-commutative structures  
with associated manifolds was studied.    

In this paper we consider chain total double complexes
 of non-abelian simplexes with 
associated spaces and reductive differentials. 
The reductivity property explained in the text allows to 
prove the relativity of corresponding cohomology, as well as its equivalence to 
 coset spaces of functionals associated to differentials of chain complexes. 
The main conjecture \ref{conjecture}
 is illustrated by the explicit proof of Theorem \ref{theorem2} 
describing a particular case 
of the simplex, the total chain double complex of 
associated spaces of correlation functions,
and intrinsic invariant bundle for a vertex operator algebra 
\cite{B, DL, FLM, FHL, FLM, Huang, K, LL}  
considered on Riemann surfaces \cite{FK}  
of various genus.  
The geometrical meaning of the theorem provides a vertex operator algebra description of 
Bott-Segal relation \cite{BS} for Lie algebras. 
\subsection{Double complex families with reductive differentials}
\label{simplex}
Let $X$ be (non necessary commutative) space of simplexes
 endowed with a double filtration 
$X=\bigcup_{\kappa, n \ge 0} X_{\kappa, n}$,
with an associated functional space $\mathcal C^{\kappa, n}(X_{\kappa, n})$. 
Let us define reductive differentials 
$\mathcal D^\kappa=\mathcal D^\kappa(X)$, $\mathcal D^n=\mathcal D^n(X)$ such that  
\begin{eqnarray}
\left(X_{\kappa+1, n}, \mathcal C^{\kappa+1, n}\right)  
= \mathcal D^\kappa. \left(X_{\kappa, n}, \mathcal C^{\kappa, n}\right),   
\quad  
\left(X_{\kappa, n+1}, \mathcal C^{\kappa, n+1}\right) 
= \mathcal D^n. \left(X_{\kappa, n}, \mathcal C^{\kappa, n}\right).  
\end{eqnarray}
Requiring single chain complex property 
with respect to each of the differentials 
\begin{eqnarray}
 \mathcal D^{\kappa+1} \circ \mathcal D^\kappa.\left(\mathcal C^{\kappa, n}\right)=0, 
\qquad 
\mathcal D^{\n+1} \circ \mathcal D^n.\left(\mathcal C^{\kappa, n}\right)=0, 
\end{eqnarray}
and the double complex property 
\begin{equation}
\label{surok}
\left(\mathcal D^\kappa \circ\mathcal D^n 
-\mathcal D^n \circ \mathcal D^\kappa\right).\left(\mathcal C^{\kappa, n}\right)=0, 
\end{equation}
the diagram 
\begin{eqnarray}
\label{buzovaish}
&& \qquad \quad \vdots \qquad \qquad \qquad \vdots 
\nn
&& \qquad \quad \downarrow \mathcal D^{\kappa-1} \qquad \downarrow \mathcal D^{\kappa-1} 
\nn
&& \cdots  \rmap  \mathcal C^{\kappa, n}
  \stackrel{\mathcal D^{\kappa, n}} {\rmap}  \mathcal C^{\kappa, n+1}  
  \rmap
 \cdots  
\nn
&&  \qquad \quad \downarrow \mathcal D^\kappa \qquad \qquad \downarrow \mathcal D^\kappa 
\nn
&& \cdots \rmap  \mathcal C^{\kappa+1, n}  \stackrel{\mathcal D^{\kappa+1, n}}  {\rmap}  
\mathcal C^{\kappa+1, n+1}
 \to \cdots  
\nn
&&  \qquad \quad \downarrow \mathcal D^{\kappa+1} \qquad \downarrow \mathcal D^{\kappa+1} 
\nn
&& \qquad \quad \vdots \qquad \qquad \qquad \vdots
\end{eqnarray} 
is then commutative.  
For countable direct sums of $C^{\kappa, n}(X_{\kappa, n})$,  
 one introduces the total complex 
$\mathcal C^m=\bigoplus_{m=\kappa+n} \mathcal C^{\kappa, n}(X_{\kappa, n})$,   
Corresponding differential is given by 
\begin{eqnarray}
d^m=\mathcal D^\kappa +(-1)^\kappa \mathcal D^n, 
\qquad 
d^m\circ d^{m-1}.\left(\mathcal C^m\right)=0, 
\end{eqnarray}
with the cohomology of the total complex $(d^m, \mathcal C^m)$ defined   
in the standard way.  
We call the single part of a complex $\mathcal C^{\kappa, n}$ 
reductive if 
 $\mathcal C^k= \mathcal D\circ \cdots \circ \mathcal D. \mathcal C^0
=P(k).\mathcal C^0$
with some operators $P(k)$, $k \ge 0$, 
and $\mathcal D$ is a finite combination 
of $\mathcal D^\kappa$ and $\mathcal D^n$. 

Introduce now the maps $\Phi$, $\Psi: X\to Y$;  
 $F$, $\mathcal G: Y \to W$, 
 the action $\Psi.\Phi: Y\times Y \to Y$,   
for spaces $Y$, $W$, and  
 a map $G: Y \to W\times W$, 
 $x\in X$, of the form 
\begin{eqnarray}
\label{partos0}
G(\Psi, \Phi)= F(\Psi(x')).\mathcal G(\Phi(x)) 
 + F(\Phi(x)).\mathcal G(\Psi(x')) +  
\sum_{{x'_0, x_0} \atop {\subset X}} \mathcal G \left(F(\Psi(x'_0)).\Phi(x_0) \right). &&    
\end{eqnarray}    
For a double-filtered $X$ 
 denote by 
$G^m=\left\{\bigoplus_{m=\kappa+n} G(x), x\in X_{\kappa, n} \right\}$,   
 the space of functionals $G(\Psi,\Phi)$ satisfying \ref{partos0}, and  
 by $Con^m$ the space of vanishing $G(\Psi,\Phi)$.     
 Let us fix the maps $\mathcal G$, $\Phi$, $\Psi$. 
  $G(\Psi, \Phi)$ depends on the map $F$ as a functional. 
Suppose that the differential of the total complex 
$d=d(F)$ is also a functional on $F$. 
If we fix a subspace $\mathfrak F$ of maps $F$, then   
the cohomology $H^m(C^m, \mathfrak F)$ of the total complex is 
a relative cohomology with respect to $\mathfrak F$. 

We call $G(\Psi, \Phi)$ covariant with respect to the differentials 
 if \eqref{partos0} 
remains of the same general form 
under arbitrary combinations of $\mathcal D^\kappa$ and $\mathcal D^n$.   
We then formulate the following conjecture which is a counterpart of 
 a proposition of \cite{BS, Wag}, i.e., the Bott-Segal theorem. 
\begin{conjecture}
\label{conjecture}
The relative cohomology of the reductive 
 chain total complex $(d^m(F)$, $C^m)$  
is equivalent to the coset space $Con^m/G^m$ 
 for some $m$. 
For $G$ covariant with respect to $d^m$,   
the equivalence extends to all $m\ge 0$.  
\end{conjecture}
Note that the vanishing \eqref{partos0} 
 represents a version of Leibniz rule. 
Thus, the cohomology relation above measures the 
inclination of $F$ and $\mathcal G$ from that rule.  

Our main example of the construction above is provided by the space $X$ of 
 $n$-simplexes of pairs $x=(v, z)$ of a vertex operator algebra elements $v\in V$
and a formal parameter $z$, and the space $\mathcal C=C^{\kappa, n}(V)$ 
  of vertex operator algebra $V$-module $W$ $n$-point  
correlation functions 
\cite{BPZ, Kn, BZF, DGM, EO, FMS, FS, TK, TUY, Zhu}. 
considered on a Riemann surface of genus $\kappa=g$. 
Due to the structure of correlations functions and reduction relations 
\cite{Zhu, BKT, MT, MT1, MT2, MT3, MT31, MTZ, Miy, T2, TZ, TZ1, TZ2, TZ3, GT, TW} 
one can form chain complexes of converging $n$-point functions. 
In this paper we assume that all $n$-point correlation functions 
are reductive to corresponding zero-point correlation function 
at any genus of Riemann surfaces. 

Recall the notion of a vertex operator algebra 
bundle given in Appendix \ref{vertex}. 
Motivated by the definition of a holomorphic connection  
  for a vertex operator algebra bundle (cf. Section 6, \cite{BZF} and \cite{Gu}) over   
a  smooth complex curve, we introduce the definition of 
the multiple point connection of the vertex operator algebra bundle (see also \cite{TUY})
over a direct product $\bigoplus_{g\ge 0} \Sigma^{(g)}$.  
With $G(\Psi, \Phi)=0$, the map $\mathcal G$ provides a generalization  
of the classical holomorphic connection over a smooth variety. 
We call the functional $G(\Psi, \Phi)$ \eqref{partos0} 
 the form of connection.  
The main results of this paper is the following theorem. 
\begin{theorem}
\label{theorem2}
The relative cohomology $H^m(W, \mathfrak F)$
 of the chain total $D^n$-reductive complex  
$(d^m$ $(F)$ , $C^m)$ of a vertex operator algebra $V$-module 
$W$ correlation functions on 
the direct product of Riemann surfaces 
is isomorphic to the factor space of 
 $D^g$- and $D^n$-covariant    
 connections $\mathcal G$    
over the space of $G^{m-1}$ of $(m-1)$-forms 
 on corresponding $V$-module $W$-bundle $\W$. 
\end{theorem}
The plan of the paper is the following.  
In Section \ref{chain} we recall the notion of vertex operator algebra $V$-module 
correlation functions, construct single chain complexes, 
and reductive differentials.  
In Section \ref{total} 
we describe the total chain complex. 
Section \ref{proof} contains a proof of the main result of this paper, 
Theorem \ref{theorem2}. 
In Appendix \ref{vertex} we recall the notion of a vertex 
operator algebra and 
its properties.  
Appendix \ref{rho} is devoted to a formalism of composing a genus $g+1$ Riemann 
surface starting from a genus $g$ Riemann surface.  
Appendix \ref{elliptic} reviews classical and generalized elliptic functions. 
In Appendix \ref{examples} examples of spaces of correlation functions 
and reduction formulas are provided. 

The results of this paper may be interesting in various fields of mathematics 
including mathematical physics \cite{DGM, FS, TK},  
Riemann surface theory \cite{FK, Gu, Kaw}, theta-functions \cite{Fa, Mu}, 
 cosimplisial geometry of manifolds \cite{Fei, Wag, BG},  
 non-commutative geometry, modular forms \cite{Miy, Be2, La, Se},  
 and the theory of foliations \cite{BG, BGG}.    
\section{The families of chain complexes for vertex operator algebra correlation functions} 
\label{chain}
\subsection{Spaces of correlation functions}
In this Section we introduce a family of chain complexes 
of correlation functions  
 a vertex operator algebra $V$-module $W$   
on a genus $g$ Riemann surface.  
Let us fix a vertex operator algebra $V$. 
Depending on its commutation relations and configuration of a genus $g$ 
Riemann surface $\Sigma^{(g)}$, the space of all $V$-module $W$ multipoint functions  
may represent various forms of complex functions defined on $\Sigma^{(g)}$.   
Consider by ${\bf v}_{n, g}=(v_{1, g}, \ldots, v_{n, g}) \in V^{\otimes n}$ 
a tuple of vertex operator algebra elements.  
Pick $n$ points on a Riemann surface $\Sigma^{(g)}$.   
Denote by ${\bf z}_{n, g}=(z_{1, g}, \ldots, z_{n, g})$ local  
coordinates around that points. 
Let us introduce our standard notation: 
${\bf x}_{n, g}= \left({\bf v}_{n, g}, {\bf z}_{n, g} \right)$.
Note that 
we use such notations to emphasize that  
the elements ${\bf x}_{n, g}$ may be chosen different for 
different genuses. 

Let $B^{(g)} \subset {\mathcal B}^{(g)}$ be moduli parameters 
 describing $\Sigma^{(g)}$. Here ${\mathcal B}^{(g)}$ 
is the set of moduli parameters for all genus $g$ Riemann surfaces. 
In particular, ${\mathcal B}^{(g)}$ characterizes a geometrical 
way if $\Sigma^{(g)}$ was constructed in a sewing procedure \cite{Y}.  
As we mentioned in Introduction, 
for each genus an $n \ge 0$-point vertex operator algebra $V$-module 
$W$ correlation function   
$\F^{(g)}_W \left({\bf x}_{n, g}; B^{(g)}\right)$ 
on $\Sigma^{(g)}$ 
has a certain specific form.
 It depends on $g$, $B^{(g)}$, 
the way a Riemann surface $\Sigma^{(g)}$ was formed,
 the type of  conformal field theory 
 model used for definitions of multipoint functions, 
 and the type of commutation relations for $V$-elements. 
We assume that, for a fixed Riemann surface set of parameters $B^{(g)}$, 
multiple point functions are completely 
determined by all choices of 
${\bf x}_{n, g} \in V^{\otimes n}\times \left(\Sigma^{(g)}\right)^n$.  
Thus, in the $\rho$-sewing procedure described in Appendix \ref{rho},  
the reduction cohomology can be treated as depending on
the set of ${\bf x}_{n, g}$ only 
with appropriate action of endomorphisms generated by $x_{n+1, g}$. 
 
For $\Sigma^{(g)}$,  
a $V$-module $W$, and $n \ge 0$, 
${\bf x}_{n, g} \in V^{\otimes n}\times \left(\Sigma^{(g)}\right)^n$,     
  we consider the spaces of all multipoint correlation functions     
$C^{g, n}(W)$ $=$ $\left\{  
\F_W^{(g)} \right.$ $\left. \left({\bf x}_{n, g}; B^{(g)}\right)  
\right\}$.    
Note that we choose elements of ${\bf v}_n$ belong to the same $V$-module $W$. 
A construction with different $V$-modules $W_i$ will be considered elsewhere.  
 Since we fix a vertex operator algebra module $W$ and $B^{(g)}$ 
we will omit them in what follows where it is possible.  
\subsection{Single chain complexes for vertex operator algebra multipoint functions} 
\label{single}
In this subsection we recall the definition of a single chain complex with respect to 
the number of points, 
and introduce a single complex with respect to the raise of genus 
of corresponding Riemann surface. 
The differentials for corresponding complexes 
are constructed according to the previous experience \cite{MT, MTZ, TZ, GT, TW, BKT}  
in applying the reduction procedure to 
vertex operator algebra $n$-point functions. 
 For $g \ge 0$, $n \ge 0$, define       
\begin{eqnarray*}
D^n(x_{n+1, g}, g): C^{g, n} & {\rightarrow }& C^{g, n+1},       
\end{eqnarray*}
\begin{eqnarray}
\label{delta_operator}
  D^n(x_{n+1, g}, g) &=& D_1^n(x_{n+1, g}, g)  + D_2^n(x_{n+1, g}, g),      
\end{eqnarray}
\begin{equation}
\label{reduction}
\F^{(g)}_W\left( {\bf x}_{n+1, g} \right)= D^n(x_{n+1, g}, g) \; 
\F^{(g)}_W\left( {\bf x}_{n, g} \right), 
\end{equation}
with differentials 
$D_1^n(x_{n+1, g}, g)$, $D_2^n(x_{n+1, g}, g)$ given by   
\begin{eqnarray}
\label{poros}
  &&D_1^n(x_{n+1, g}, g).\F^{(g)}_W \left( {\bf x}_{n, g} \right) =  
\sum\limits_{l=1}^{l(g)} 
  f_1^{(g)} \left( x_{n+1, g}, l \right)  
\;  T^{(g)}_l(x_{n+1, g}). \F^{(g)}_W \left( {\bf x}_{n, g} \right), 
\\
 &&D_2^n (x_{n+1, g}, g). \F^{(g)}_W \left( {\bf x}_{n, g} \right) 
 =
\nonumber
 \sum\limits_{k=1}^{n} \sum\limits_{m \ge 0}   
  f_2^{(g)} ( x_{n+1, g}, k, m ) \;   
T_k^{(g)}( v_{n+1, g}(m)).\F^{(g)}_W \left( {\bf x}_{n, g} \right), 
\end{eqnarray} 
where $l(g) \ge 0$ is a constant depending on $g$, 
and the meaning of indexes $1 \le k \le n$, 
$1 \le l \le l(g)$, $m \ge 0$ explained below. 
Then the operator 
$T^{(g)}_l(x_{n+1}).\F^{(g)}_W \left( {\bf x}_{n, g} \right)$  
gives a function of $F^{(g)} \left( {\bf x}_{n, g} \right)$ 
depending on $x_{n+1, g}$. 
The operator 
\begin{equation*}
T^{(g)}_k(v_{n+1, g}(m)).\F^{(g)}_W \left( {\bf x}_{n, g}  \right) =   
\F^{(g)}_W \left( T_k(v_{n+1, g}(m)).{\bf x}_{n, g} \right), 
\end{equation*}
is the insertion of  
the $m$-th mode $v_{n+1, g}(m)$ (or $v_{n+1, g}[m]$-mode depending on $g$),
 $m \ge 0$. 
of vertex operator algebra elements $v_{n+1, g}$, 
 in front of the $k$-th argument $v_{k, g}$ of $x_{k, g}$ 
inside the $k$-th vertex operator in the functional 
 $\F^{(g)}_W \left( {\bf x}_{n, g} \right)$.  
Here we use the notation 
\[
T^{(g)}_k(\gamma).\; f({\bf x}_{n, g}) =   
f\left(x_{1, g},  \ldots,  \gamma.x_{k, g},  \ldots, x_{n, g} \right), 
\]
for an operator $\gamma$ acting on $k$-th argument of a functions $f$. 
 Note that commutation properties of $D_1^n(x_{n+1, g}, g)$ 
 and $D_2^n(x_{n+1, g}, g)$ depend on genus $g$.  
 Operator-valued functions $f^{(g)}_1\left( x_{n+1, g}, l \right) \; 
T^{(g)}_l( v_{n+1, g})$,  
$f_2^{(g)}(x_{n+1, g}, k, m)$. $T^{(g)}_k( v_{n+1, g}(m))$    
 depend on genus of a Riemann surface $\Sigma^{(g)}$.  
For $n \ge 0$, 
let us denote by ${\mathfrak V}_n$  
the subsets of all
 $x_{n+1, g} \in V\times \Sigma^{(g)}$, 
 such that  
 the chain condition 
\begin{eqnarray}
\label{nconditions}
D^{n+1}(x_{n+2, g}, g)\circ 
D^n(x_{n+1, g}, g). \F^{(g)}_W \left({\bf x}_{n, g} \right)=0,   
\end{eqnarray} 
 for the differentials \eqref{poros} for complexes $C^{g, n}$   
 is satisfied.    
 
Next, consider the differentials  
\begin{eqnarray}
\label{operabase}
D^g: C^{g, n} \to C^{g+1, n}, \quad 
 \F^{(g+1)}_W\left({\bf x}_{n, g+1}; B^{(g+1)}\right)  
= D^g.\F^{(g)}_W \left({\bf x}_{n, g}; B^{(g)}\right), &&    
\end{eqnarray}
for $V$-module $W$ on a genus $g$ Riemann surface. 
There exist \cite{Y, MT, MT2, MT3, MT31, TZ, TZ1, GT, TW, Be2, Kn}  
various geometrical ways how to 
increase the genus of a Riemann surface, and, therefore,  
 ways how to introduce corresponding differential $D^g$. 
In this paper we will use the $\rho$-formalism of attaching a handle to a genus $g$
Riemann surface to form a genus $g+1$ Riemann surface (see Appendix \ref{rho}). 
In this geometric setup the differential $D^g$ is given by 
\begin{eqnarray}
\label{rhoraise}
&&\F^{(g+1)}_W \left({\bf x}_{n, g+1}; B^{(g+1)}\right)
= D^g.\F^{(g)}_W \left({\bf x}_{n, g}; B^{(g)}\right), \qquad 
\nn
 &&\F^{(g+1)}_W \left({\bf x}_{n, g+1}; B^{(g+1)}\right) 
= \sum_{k \ge 0} \; \sum\limits_{w_k \in W_{(k)}} \rho_g^k \;   
  \F^{(g)}_W \left({\bf x}_{n, g}, \overline{w}_k,
 \zeta_1, w_k, \zeta_2; B^{(g)}\right) 
\nn
&&\qquad \qquad \qquad \qquad \quad 
=  \sum_{k \ge 0} \sum\limits_{w_k \in W_{(k)}} \rho_g^k   
 T(\overline{w}_k, \zeta_1, w_k, \zeta_2).  
\F^{(g)}_W \left({\bf x}_{n, g}; B^{(g)}\right), 
\end{eqnarray}
Note that in this formulation the differential $D^g$ does not depend on $n$. 
It is assumed that \ref{rhoraise} converges in $\rho$ for $W$.  
The resulting expression for $\F^{(g+1)}_W$ $\left({\bf x}_{n, g+1}\right)$
depends on the positions of 
$(\overline{w}_k$, $\zeta_1)$, $(w_k$, $\zeta_2)$-insertions  
into $\F^{(g)}_W$ $ \left({\bf x}_{n, g}\right)$ and their permutation properties 
with ${\bf x}_{n, g}$. 
Here we fix the position of insertion right agter the element $x_{n, g}$
 as it was done in \cite{TZ, T2}.   

In this paper we consider $C^{g, n}$ as spaces of arbitrary 
 $\F^{(g)}_W\left({\bf x}_{n, g} \right)$ not necessary obtained 
as a result of $\rho$-procedure from some 
$\F^{(g-1)}_W\left({\bf x}_{n', g-1}\right)$ considered on a 
genus $g-1$ Riemann surface. 
At the same time we act on 
$\F^{(g)}_W\left({\bf x}_{n, g} \right)$ by the differential $D^g$
which involves the $\rho_g$-sewing procedure.
We assume also that $\F^{(g)}_W\left({\bf x}_{n+1, g}\right)$ 
can be obtained from $\F^{(g)}_W\left({\bf x}_{n, g} \right)$
via reduction formulas. 
The single chain complex condition
 for the differentials $D^g$ \eqref{rhoraise} has the form 
\begin{equation}
\label{gconditions}
D^{g+1}\circ D^g. \F^{(g)}_W\left({\bf x}_{n, g} \right)=0. 
\end{equation}
For $g \ge 0$,  
let us denote by ${\mathfrak V}_g$  
the subsets of all
 $x_{n, g} \in V \times \Sigma^{(g)}$,  
 such that  
 the single chain condition \eqref{gconditions}  
 for the differentials \eqref{rhoraise} of complexes $C^{g, n}$    
 is satisfied.    

For the single chain complexes with differentials 
$D^g$ and $D^n(x_{n+1, g}, g)$ 
 given in this Section one  
defines corresponding partial cohomology in the standard way.  
Combining the formulations for two single complexes, we introduce the family 
of the single chain complexes $C^{g, n}$, $g \ge 0$, $n\ge 0$ 
for a vertex operator algebra $V$ on Riemann surfaces. 
Using the differentials \eqref{poros} and \eqref{operabase}
 one can compose the functional $G(\Psi, \Phi)$ \eqref{partos0}. 
For each $g$ and various types of vertex operator algebras, 
there exists standard sets of operators $F$ 
\cite{BKT, GT, MT, MT1, MT2, MT3, MT31, MTZ, T2, TZ, TZ1, TZ2, TZ3, TW, Zhu},  
i.e., $T(\overline{w}_k, \zeta_1, w_k, \zeta_2)$, 
$f_1^{(g)} \left( x_{n+1, g}, l \right) \;T^{(g)}_l(x_{n+1, g})$, 
and $f_2^{(g)} ( x_{n+1, g}, k, m)$  
 $T_k^{(g)}( v_{n+1, g}(m)).$  
determining the differentials 
$D^g$, $D^n(x_{n+1, g}, g)$.   
\section{The chain total complex}
\label{total}
In order to turn the family $C^{g, n}$ of single chain complexes 
into a double chain complex we have to apply further requirement of commutation on 
 the differentials $D^g$ and $D^n(x_{n+1, g}, g)$.  
In addition to the conditions \eqref{nconditions} and \eqref{gconditions}
we require for the differentials $D^g$ and $D^n(x_{n+1, g}, g)$ to satisfy 
\begin{eqnarray}
\label{gnconditions}
&& \left(D^g \circ  D^n(x_{n+1, g}, g) 
- D^n(x_{n+1, g}, g) \circ D^g  \right).   
\F^{(g)}_W\left({\bf x}_{n, g} \right)=0.  
\end{eqnarray}
Then the family $C^{g, n}$ turns into a chain double complex. 
We denote by ${\mathfrak V}_{g, n}$  
the subsets of all $x_{n+1, g}  
\in V\times \Sigma^{(g)}$,  
 such that \eqref{gnconditions} is satisfied. 
\subsection{The total complex}
\label{directtotal}
For the double complex $C^{g,n}$ 
the associated total complex 
 is given by
${\rm Tot}^m (C^{g,n}) = \bigoplus_{m=g+n} C^{g,n}=C^m$,    
for $m \ge 0$, with the differential 
\begin{eqnarray}
\label{totalcomp}
 d^m: {\rm Tot}^m (C^{g,n}) \to {\rm Tot}^{m+1} (C^{g,n}),  \quad  
 d^m = \sum_{m=g+n} \left(D^g+(-1)^g D^n (x_{n+1, g}, g)\right). &&   
\end{eqnarray}
Note that the differential $D^n(x_{n+1,g}, g)$ 
in \eqref{totalcomp} is defined for all choices of
 $x_{n+1, g} \in V \times \Sigma^{(g)}$,  
  chosen separately for all possible combinations of $g$ and $n$ 
such that $m=g+n$. 
For $n \ge 0$, $g \ge 0$, $m=g+n$, 
let us denote by ${\mathfrak V}_d$ 
the subsets of all
 $x_{n+1, g} 
\in V \times \Sigma^{(g)}$,  
 such that  
 the chain conditions \eqref{nconditions},
 \eqref{gconditions}, \eqref{gnconditions}  
and 
\begin{equation}
\label{totalcondition}
d^{m+1}\circ d^{m}.\left(C^m\right)=0, 
\end{equation} 
 for the differentials \eqref{poros} for complexes $C^{g, n}$   
 are satisfied.    
Note that ${\mathfrak V}_g$, 
${\mathfrak V}_n$, and ${\mathfrak V}_{g, n}$ are subsets of ${\mathfrak V}_d$. 
The spaces with conditions \eqref{nconditions}, 
\eqref{gconditions}, \eqref{gnconditions}, and \eqref{totalcondition}
constitute a semi-infinite chain 
double complex with the commutative diagram 
\begin{eqnarray*}
 &&\quad \qquad \vdots \qquad \qquad \qquad \qquad \quad \vdots 
\nn
 &&\qquad \quad \downarrow D^{g-1} \qquad \qquad \qquad \downarrow  D^{g-1} \ldots 
\nn
 &&0  \rmap  C^{g, 0}  \qquad \stackrel{D^0(x_{1, g}, g)}  {\rmap} \qquad C^{g, 1} 
 \stackrel{D^1(x_{2, g}, g)}  {\rmap} 
 \ldots  \rmap  C^{g, n-1}  \stackrel{D^{n-1}(x_{n, g}, g)}  {\rmap}  C^{g, n} 
 \rmap   
\nn
  &&\qquad \quad \downarrow D^g \qquad \qquad \qquad \quad \downarrow D^g \qquad \ldots 
\nn
 &&0  \rmap  C^{g+1, 0}  \stackrel{D^0(x_{1, g+1}, g+1)}  {\rmap}  C^{g+1, 1}
 \stackrel{D^1(x_{2, g+1}, g+1)}  {\rmap} 
 \ldots  
\stackrel{D^{n-1}(x_{n, g+1}, g+1)} {\rmap} C^{g+1, n} 
 \rmap   
\nn
  &&\qquad \quad \downarrow D^{g+1} \qquad \qquad \qquad  \downarrow D^{g+1}  \qquad \ldots 
\nn
 &&\quad \qquad \vdots \qquad \qquad \qquad \qquad \quad \vdots 
\end{eqnarray*} 
and standardly defined the reduction cohomology involving $d^m$  
of the total complex \eqref{totalcomp}. 

Due to vertex operator algebra properties, 
the conditions \eqref{nconditions}, \eqref{gconditions}, \eqref{gnconditions}, and 
  \eqref{totalcondition} result in expressions
 containing finite series of vertex operator algebra modes 
and coefficient functions. That conditions 
 narrow the space of compatible elements $x_{n+1, g}$, and, therefore, 
 corresponding multipoint functions $\F^{(g)}({\bf x}_{n, g})$.   
Nevertheless, the subspaces of $C^{g, n}(W)$, $g \ge 0$, $n \ge 0$, 
of multipoint functions 
such that the conditions above are fulfilled  
for reduction cohomology complexes 
are non-empty. 
 For all $g$,  
 the conditions mentioned represent  
an infinite $n \ge 0$, $g \ge 0$ set of functional-differential  
equations (with finite number of summands) on converging complex functions 
$\F^{(g)}_W\left({\bf x}_{n, g} \right)$ defined for   
 $n$ local complex variables on Riemann surfaces 
of genus $g$ with extra action of the operator 
$T(\overline{w}_k, \zeta_1, w_k, \zeta_2)$, 
with functional coefficients
$f_1^{(g)} \left( x_{n+1, g}, l \right)$,  
 $f_2^{(g)} ( x_{n+1, g}, k, m)$.    
In examples given in Appendix \ref{examples}  
the functional coefficients
$f_1^{(g)} \left( x_{n+1, g}, l \right)$,  
 $f_2^{(g)} ( x_{n+1, g}, k, m)$
 are genus $g$ generalizations of elliptic 
 functions
 on $\Sigma^{(g)}$. 
Note that 
all vertex operator algebra elements of ${\bf v}_n\in V^{\otimes n}$,   
 as non-commutative parameters  
are not present in final form of functional-differential equations 
since they incorporated into either matrix elements, traces, and other 
forms in corresponding genus $g$ multipoint functions.  
According to the theory of such equations, 
equations resulting from 
\eqref{nconditions} \eqref{gconditions}, 
\eqref{gnconditions}, and  \eqref{totalcondition} 
always have non-vanishing solutions in the domains 
they are defined. 
Applying the reduction procedure by differentials 
$D^n(x_{n+1, g}, g)$ we reduce the functions 
$\F^{(g)}_W\left({\bf x}_{n, g}\right)$ 
to corresponding zero-point functions $\F^{(g)}_{W, 0}$, 
i.e., we obtain in general 
$\F^{(g)}_W\left({\bf x}_{n, g}\right)
= P_n\left(F; g, n,{\bf x}_{n, g}\right)\;\F^{(g)}_{W, 0}$, 
where $P_n\left(F; g, n, {\bf x}_{n, g}\right)$ 
are explicitly computable functions 
 containing genus $g \ge 0$ 
generalized elliptic functions (see Appendix \ref{elliptic}). 
For non-zero zero-point functions, equations
\eqref{nconditions}, \eqref{gconditions}, \eqref{gnconditions}, and 
  \eqref{totalcondition}
expressed for $P_n\left(F; g, n,{\bf x}_{n, g}\right)$  
 can be solved by methods of the analytic number theory. 
\section{Proof of the main theorem}    
\label{proof}
In this Section we provide a proof of Theorem \ref{theorem2}.  
First, recall that the form of $\F^{(g)}_W\left({\bf x}_{n, g} \right)$  
is specific for each $g$. 
The differential $D^g$ \eqref{operabase} 
makes the genus transitions from $\F^{(g)}_W\left({\bf x}_{n, g}\right)$ 
to $\F^{(g+1)}_W\left({\bf x}_{n, g+1}\right)$.    
Recall the notion of a vertex operator algebra bundle given in Appendix \ref{vertex}. 
The definition \eqref{poetika} corresponds to the case $g=0$ of the $\W^*$-section of 
the vertex operator algebra $V$ bundle.  
One can see that \eqref{poetika} combined with the differential $D^g$ \eqref{operabase}  
extends \eqref{poetika} to $g>0$ cases.   
Namely, we define $\F^{(g)}_{\W, n}({\bf x}_n)$ as follows 
\begin{eqnarray}
\label{baza}
 \F^{(g)}_\W({\bf x}_n)&=& \sum\limits_{ {\bf k}_g \in \Z^g, \; u_k \in V_{(k)}  }  
\langle ( \overline{u}_k, {\bf z}_n; \overline{w}_{ k_j }, \zeta_{1, j}),  
\mathcal Y_\W( {\bf v}_n)\cdot (u_k, {\bf z}_n; w_{  k_j } , \zeta_{2, j})\rangle     
\nn
&=& D^g \circ \cdots \circ D^0. \F^{(0)}_\W({\bf x}_n), 
\end{eqnarray}
for $\zeta_{a, j}$, $a=1$, $2$, 
 $u_k \in W_{(k)}$, $w_{k_j} \in W_{(k_j)}$, $1 \le j \le g$, 
  $\overline{u}_k$, $\overline{w}_{k_j}$ their corresponding duals,  
and where the action of $D^g$, $\ldots$, $D^0$ is realized 
via \eqref{operabase}.   
Note that \eqref{baza} preserves the linearity properties 
 in $\overline{u}$, $u$, and $\mathcal O_x$-linearity in ${\bf v}_n$.  
Let us mention that due to the relation \eqref{baza} we could formulate 
all the material of Sections \ref{chain}--\ref{total} 
in terms of complexes constituted by 
$\F_{\W, m}$. 

Let us denote by $\mathcal{SW}$ the space of sections of 
the vertex operator algebra $V$-module $W$ bundle $\W$.  
In our setup, 
 by identifying $\Psi$ and $\Phi$ with 
sections $\psi(x')$ and $\phi(x)$ of $\W$ correspondingly,   
the map $G(\Psi, \Phi)$ \eqref{partos0} 
with a $\C$-multi-linear map 
$\mathcal G: \mathcal{SW}= \bigoplus_{m=g+n} \W^{\otimes n}   
  \times \left(\Sigma^{(g)} \right)^n  \to \C$,  
and for any operator $F$, 
turns into 
\begin{eqnarray}
\label{locus}
G(\psi({\bf x}'), \phi({\bf x}))=  F(\psi({\bf x}')). 
\mathcal G \left(\phi({\bf x}) \right) + F(\phi({\bf x})). 
\mathcal G\left(\psi({\bf x}') \right)    
+ \mathcal G\left( F(\psi({\bf x}')).\phi({\bf x}) \right). 
\end{eqnarray}
The vanishing map \eqref{locus}  
gives raise to 
 a generalization $\mathcal G$ of the holomorphic connection  
on $\W$. 
Geometrically, for a vector bundle $\W$ defined over $\Sigma^{(g)}$,  
the vanishing generalized connection \eqref{locus} relates 
two sections $\psi({\bf x}')$ and $\phi({\bf x})$.  
Now we are ready to give a proof of Theorem \ref{theorem2} 
\begin{proof} 
Let us denote $\F_{\W, m}=\sum_{m=g+n} 
\F^{(g)}_\W\left({\bf x}_{n, g} \right)$.    
We assume that operators $F$ satisfy 
\eqref{nconditions}, \eqref{gconditions}, \eqref{gnconditions}, 
and \eqref{totalcondition}.  
The definition \eqref{baza} provides the coordinateless 
expression for $\F_{\W, m}$,  
i.e., 
\begin{eqnarray}
\label{poetikanew}
\F_{\W, m}= \sum\limits_{m=g+n} \; \sum\limits_{ {\bf k}_g \in \Z^g}  
\sum\limits_{u_k\in V_{(k)}} 
\langle  ( \overline{u}_k, {\bf z}_n; \overline{w}_{ k_j}, \zeta_{1, j}),  
 \Y_\W(i_z ({\bf v}_n)) \cdot 
(u_k, {\bf z}_n; w_{ k_j}, \zeta_{2, j}, ) \rangle. &&  
\end{eqnarray}
Using \eqref{poros} and  \eqref{rhoraise} we set for $m=g+n$,  
$m+1=g'+n'$, 
\begin{eqnarray*}
\mathcal G\left(\phi({\bf x})\right) =\F_{\W, m}, \;   
 \psi({\bf x}')= \Y_\W(i_z ({\bf v}_{n', g'})) \cdot (., {\bf z}_{n', g'}), \;     
 \phi({\bf x})= \Y_\W(i_z ({\bf v}_{n, g})) \cdot (., {\bf z}_{n, g}), && 
\end{eqnarray*}
\begin{eqnarray}
\label{identificationsforg} 
 -  F\left(\psi({\bf x}') \right).  
\mathcal G \left(\phi({\bf x}) \right)=
 \left( \sum_{k \ge 1} \; \sum\limits_{w_k \in W_{(k)}} \rho_g^k \;  
 T\left(\overline{w}_k, \zeta_1, w_k, \zeta_2\right)\right).\F_{\W, m}, 
\qquad \qquad \qquad  && 
\\
  -\mathcal G\left( F(\psi({\bf x}')).\phi({\bf x}) \right)  
=(-1)^g\left[\sum\limits_{l=1}^{l(g)} 
  f^{(g)}_1\left(x_{n+1, g}, l \right)  \;  T^{(g)}_l(x_{n+1, g})
 \right. \qquad \qquad \qquad  &&
\nn
  \left. \qquad \qquad \qquad \qquad 
+\sum\limits_{k=1}^n \sum\limits_{r \ge 0}   
  f_2^{(g)} ( x_{n+1, g}, k, r) \; 
T_k^{(g)}( v_{n+1, g}(r)) \right].\F_{\W, m}. &&
\nonumber
\end{eqnarray}
Now, let us assume that the form 
of the generalized connection $\mathcal G=\F_{\W, m}$ 
remains the same 
for the vanishing $G^m \in Con^m$ and non-vanishing  
  $G^m \in G^m$ given by \eqref{locus},  but  
 the operator $F$,
 and the differentials $\left(D^n(x_{n+1, g}, g)\right)'$  
 may differ from \eqref{identificationsforg}, 
 \eqref{poros}--\eqref{rhoraise}, and satisfying conditions 
\eqref{gconditions}, \eqref{nconditions}, \eqref{gnconditions},   
and \eqref{totalcondition}. 
By using the reduction procedure given by some other operators $F'$ 
and some different differentials $(D^n(x_{n+1, g}, g))'$,     
we reduce \eqref{locus} to 
\begin{equation}
\label{pusto}
G(\psi({\bf x}'), \phi({\bf x}))
= \sum\limits_{m=g+n} P_m\left(F'; g, n, {\bf x}_{n, g}\right)\; \F^{(g)}_{\W, 0}, 
\end{equation}
where $P_m\left(F'; g, n, {\bf x}_{n, g}\right)$ 
are functions depending on operators $F'$, and explicitly  
 containing genus $g \ge 0$ 
generalized elliptic functions. 
Thus, $\F^{(g)}_{\W, n}$ is explicitly known and it  
is represented as a series of auxiliary functions 
$P_m\left(F'; g, n, {\bf x}_{n, g}\right)$   
depending on $F'$, $g$, $n$, and ${\bf x}_{n, g}$. 

Now let us consequently apply the initial 
differential $D^n(x_{n+1, g}, g)$ \eqref{poros}  
 to reconstruct back some functions $\F''_\W$ starting from each of 
$\F^{(g)}_\W$ for $m=g+n$ in \eqref{pusto}.  
Finally, we obtain \eqref{locus} for $\mathcal G=\F''_{\W, m}$ 
with 
\[
\F''_{\W, m}= \sum\limits_{m=g+n} D^{n-1}(x_{n-1, g})\circ \cdots \circ D^1(x_{1, g}) 
\circ P_m(F'; g, n, {\bf x}_{n, g}). \F^{(g)}_{\W, 0}.    
\]
Since the differentials $D^g$ and $D^n(x_{n+1, g}, g)$ 
are defined in that way they act on the functions 
$\F^{(g)}_{\W}$ only, it is clear that the action of 
$D^{n-1}(x_{n-1, g})\circ \ldots \circ D^1(x_1, g)$ 
and multiplication by $P_m(F'; g, n, {\bf x}_{n, g})$ commute.  

Thus, we infer that 
the $m$-th reduction cohomology $H^m$ of $(d^m, C^m)$ is 
  equivalent to 
the factor space $Com^m/G^m$ with the coefficients given by the coset 
\[
\left\{
P_m(F'; g, n, {\bf x}_{n, g})|_{G(x', x)=0}
/P_{m-1}(F'; g, n, {\bf x}_{n, g})|_{G(x', x)\ne 0}\right\},  
\]  
of genus $g$ counterparts of 
elliptic functions, and relative to the 
 subspace of covariance preserving operators $F'$ for 
the the functional $G\left(\psi({\bf x}'), \phi({\bf x})\right)$.   
The $H^m$-th relative reduction cohomology 
of a vertex operator algebra $V$-module $W$ 
is then given by the ratio of series of generalized elliptic functions 
 recursively generated by the reduction formulas 
\eqref{reduction}--\eqref{operabase}. 
\end{proof}
\section*{Acknowledgment}
The author is supported by the Academy of Sciences of the Czech Republic (RVO 67985840). 
\section{Appendix: Vertex operator algebras} 
\label{vertex}
In this Subsection we recall the notion of a 
vertex operator algebra $(V,Y,\mathbf{1}_V,\omega)$
\cite{B, DL, FHL, FLM, K, LL}.   
Here $V$ is a linear space  
endowed with a $\mathbb Z$-grading 
$V=\bigoplus_{r \in {\mathbb Z}} V_r$,     
 $\dim V_r<\infty$.
  The state  $0 \ne {\mathbf 1}_V \in V_0$, is called vacuum vector, 
$\omega\in V_2$ is the conformal vector with properties described below.
The vertex operator is a linear map  
$Y: V\rightarrow \mathrm{End}(V)[[z,z^{-1}]]$,  
with formal variable $z$. 
  For any vector $v\in V$, $x=(v, z)$, we have a vertex operator   
$Y(x)=\sum_{n\in {\mathbb Z}}v(n)z^{-n-1}$.   
The linear operators (which are called modes) 
$u(n):V\rightarrow V$ satisfy creativity 
$Y(v,z){\mathbf 1}_V= v +O(z)$,  
and lower truncation 
$v(n)u=0$,  
conditions for each $u$, $v\in V$ and $ n\gg 0$. 
 The vertex operators satisfy an analogue of Jacobi identity 
\begin{eqnarray*}
 \delta(z_1, z_2, z_0) Y (x_1 )Y(x_2)    
  - 
\delta (z_2, z_1, -z_0) Y(x_2) Y(x_1)   
= z_0 z_2^{-1} \delta(z_1, z_0, z_2) 
Y \left( Y(u, z_0)v, z_2\right),&& 
\end{eqnarray*} 
for $x_1=(u, z_1)$, $x_2=(v, z_2)$, and 
$\delta(z, z', z'')=\delta((z' - z'')(z''')^{-1})$.  
These axioms imply 
locality, skew-symmetry, associativity and commutativity conditions:
\begin{eqnarray*}
(z_1-z_2)^N
Y(x_1)Y(x_2)  
&=& 
(z_1-z_2)^N 
Y(x_2)Y(x_1), 
\\
Y(u,z)v &=& e^{zL(-1)}Y(v,-z)u,
\\
(z_0+z_2)^N Y(u,z_0+z_2)Y(v,z_2)w &=& (z_0+z_2)^N Y(Y(u,z_0)v,z_2)w,  
\\
u(k)Y(v,z)- Y(v,z)u(k) &=& \sum\limits_{j\ge 0}  \left( k \atop j \right)
Y(u(j)v,z)z^{k-j},
\end{eqnarray*}
for $u$, $v$, $w\in V$ and integers  $N \gg 0 $.   
For the conformal vector $\omega$ one has 
$Y(\omega, z)=\sum_{n\in {\mathbb Z}}L(n)z^{-n-2}$,   
where $L(n)$ satisfies the Virasoro algebra with central charge $c$ 
\begin{equation}
[\,  L(m),L(n)\, ]=(m-n)L(m+n)+\frac{c}{12}(m^3 -m)\delta_{m,-n}{\rm Id}_V, 
\label{Virasoro}
\end{equation}
where ${\rm Id}_V$ is identity operator on $V$. 
Each vertex operator satisfies the translation property  
$\partial_z Y(u,z)= Y\left(L(-1)u,z\right)$.   
The Virasoro operator $L(0)$ defined a ${\mathbb Z}$-grading with 
$L(0)u=ru$,  
 for $u\in V_r$, $r\in {\mathbb Z}$, $\wt(u)=r$.  
For $v={\mathbf 1}_V$  
one has $Y({\mathbf 1}_V, z)={\rm Id}_V$.  
Note also that modes of homogeneous states 
are graded operators on $V$, i.e., for $v \in V_k$,  
 $v(n): V_m \rightarrow V_{m+k-n-1}$. 
 In particular, let us define the zero mode $o(v)$ 
as $o(v) = v(wt (v) - 1)$ additively extending to $V$. 

 Let us recall also the square-bracket formalism \cite{Zhu} 
for a vertex operator algebra $V$, i.e., the quadruple 
$(V,Y[.,.],\mathbf{1}_V, \tilde{\omega})$ 
 The square bracket vertex operators are given by 
\[
Y[v,z]=\sum_{n\in \mathbb{Z}}v[n]z^{-n-1}=Y(q_z^{L(0)}v,q_z-1),
\]
with $q_z=e^z$. Corresponding conformal vector is 
$\tilde{\omega} =\omega -\frac{c_V}{24}\mathbf{1}_V$.  
For $v$ of $L(0)$ weight $wt(v)\in \mathbb{R}$ and $m\geq 0$, 
\[
v[m] 
=
m!\sum\limits_{i\geq m}c(wt(v),i,m)v(i),\qquad 
\sum\limits_{m=0}^{i}c(wt(v),i,m)x^{m} 
=
  \left( wt(v)-1+x \atop i \right).
\]

Given a vertex operator algebra $V$, one defines the 
 adjoint vertex operator with respect to $\rho \in \C$,    
$ Y_\rho^\dagger[v,z] = Y\left[\exp\left(z\rho^{-1} L[1]\right)
\left(- \rho z^{-2} \right)^{L[0]}v,   
 \rho z^{-1}\right]$.
associated with the formal M\"obius map \cite{FHL}
$z \mapsto \frac{\rho}{z}$.    
An element $u\in V$ is called quasiprimary if 
$L(1)u=0$.  
For quasiprimary $u$  
of weight $\wt(u)$ one has
$u^\dagger(n) = (-1)^{\wt(u)}$ $\alpha^{n+1-\wt(u)}$ $u(2\wt(u)-n-2)$. 

We call a bilinear form 
$\langle . , . \rangle : V \times V\rightarrow \C$,  
  invariant if \cite{FHL, Li} 
$\langle Y(u,z)a, b \rangle = \langle a, Y^\dagger(u,z)b \rangle$,  
for all  
$a$, $b$, $u\in V$.
Note that the adjoint vertex operator 
$Y^\dagger(.,.)$ as well as the bilinear form $\langle. , . \rangle$,   
 depend on $\alpha$. 
 Rewriting  in terms of modes, we obtain 
$\langle u(n)a,b \rangle = \langle a,u^\dagger(n)b \rangle$. 
Choosing $u=\omega$, and for $n=1$ implies that 
$\langle L(0)a,b \rangle = \langle a,L(0)b \rangle$.
Thus, $\langle a, b \rangle=0$, 
 when $\wt(a) \neq \wt(b)$.  
A vertex operator algebra $V$ is called of strong-type 
if $V_0 = \C\mathbf{1}_V$,   
and it is simple and self-dual, 
i.e., isomorphic to the dual module 
$V^\prime$ as a $V$-module.   
It is proven in \cite{Li} 
that a strong-type vertex operator algebra $V$ 
has a unique invariant non-degenerate  
bilinear form up to normalization. 
The form $\langle  ., .\rangle$ defined 
 on a strong-type vertex operator algebra $V$ 
is the unique invariant bilinear form $\langle . , . \rangle$   
normalized by 
$\langle \mathbf{1}_V, \mathbf{1}_V \rangle = 1$.
A vertex operator algebra $V$-module $W$ possesses similar properties as $V$ 
\cite{B, K, FHL, FLM, LL}. 

A vertex algebra $V$ is called quasi-conformal \cite{BZF} if 
  it admits an action of the local Lie algebra  
of ${\rm Aut}\; \Oo$ for which 
\[
[{\bm v}, Y (u, w)] = -\sum_{m \ge -1} ((m+1)!)^{-1} 
(\partial_w^{m+1} v(w))Y (L_m u, w), 
\]
with ${\bm v} = -\sum_{r \ge -1}v_rL_r$, 
$v(z)\partial_z =\sum\limits_{r\ge-1} v_r z^{r+1} \partial_z$,
is true for any  
$v \in  V$, the element $L_W(-1) = - \partial_{z}$, 
as the translation operator $T$, 
$L_W(0) = - z \partial_{z}$. 
In addition it 
  acts semi-simply with integral
eigenvalues, and the Lie subalgebra of the positive 
part of local Lie algebra of ${\rm Aut}\; \Oo^{(n)}$
 acts locally nilpotently. 
\subsection{Vertex operator algebra $V$-module $W$ bundle $\W$} 
\label{bundle}  
The notion of a vertex operator algebra $V$ bundle was introduced in \cite{BZF}.  
In this Appendix we recall that definition of a $V$-module $W$ bundle $\W$.  
The idea is to associate 
 canonically (i.e., coordinate independently) 
${\rm End} \; \W$-valued sections $\mathcal Y_\W$ of $\W^*$   
 (the bundle dual to $\W$)    
to matrix elements of $V$-module $W$ vertex operators. 

Denote by ${\rm Aut} \; \mathcal O$ the group of continuous automorphisms
of $\mathcal O = \C[[z]]$ on an arbitrary smooth curve $S$ and 
 its Lie algebra ${\it Aut} \; \mathcal O$.  
Let $V$ be a quasi-conformal vertex algebra (see the previous subsection).
 Its module $W$ is graded by   
 finite dimensional ${\rm Aut} \; \mathcal O$-submodules. 
One defines a vertex operator algebra bundle $\mathcal W_S$ 
 and its dual $\mathcal W^*_S$ as inductive and 
projective limits of vector bundles of finite rank 
 over $S$, 
 in particular, for the disc $D = {\rm Spec} \C[[z]]$, 
 or $D_z = {\rm Spec} \mathcal O_z$ with 
${\it Aut}_{D_z} ={\it Aut} \; X|_{D_z}$ and 
$\mathcal W_{D_z} = \mathcal W_{S|_{D_z}}$. 
Let ${\it Aut}_z$ be the ${\rm Aut}\; \mathcal O$-torsor 
of coordinates at $z \in S$. Recall that 
$\mathcal W_z = {\it Aut}_z {\times \atop {{\rm Aut} \mathcal O}} W$ 
is the fiber of $\mathcal W_{|D_z}$ at $z \in S$.  
Let us define ${\rm End} \; \mathcal W_z$-valued meromorphic
section $\mathcal Y_\W$ of the bundle $\mathcal W^*$  
on the punctured disc $D^\times_z$.  
This section is given by the map (linear in $\overline{u}$, $u$, 
and $\mathcal O_z$-linear in $s$) 
$(\overline{u}, s, u)  
\mapsto \langle \overline{u}, \mathcal Y_\W(s) \cdot u\rangle$,   
assigning a function on $D^\times_z$  
 denoted by $\langle \overline{u}, \mathcal Y_\W(s) \cdot v \rangle$, 
 for $\overline{u} \in \mathcal W^*_z$, $u \in \mathcal W_z$,  
 and a regular section $s$ of $\mathcal W|_{D_z}$. 
For a coordinate $z$ on the disc $D_z$, we then 
obtain a $z$-trivialization of $\mathcal W$ 
$i_z: W [[z]] \simeq \Gamma(D_z, \mathcal W)$,  
and trivializations $\mathcal W^* \simeq \mathcal W^*_z$,  
$\mathcal W \simeq \mathcal W_z$ of the fibers which we denote 
by $(\overline{u}, z)$,  $(u, z)$.  
Define an ${\rm End} \; \mathcal W_z$-valued 
section $\mathcal Y_\W$ of $\mathcal W^*$ 
 on $D_z^\times$ by  
\begin{eqnarray}
\label{poetika}
\F^{(0)}_\W(x)= \sum\limits_{w_k \in V_{(k)}}   
\langle (\overline{u}_k, z), \Y_\W(i_z (v)) \cdot (u_k, z)\rangle  
 \sim \sum\limits_{w_k \in W_{(k)}}\langle \overline{u}_k, Y(v, z)u_k \rangle
=\F^{(0)}_W(x), && 
\end{eqnarray}
where $z$ is a coordinate on $D_z$. Then the section $\mathcal Y_\W$  
is canonical, i.e., 
independent of the choice of coordinate $z$ on $D_z$.  

Let $S$ be a smooth complex variety and $\mathcal E \to S$ 
 a holomorphic vector bundle over $S$.
 We use the same notation 
$\mathcal E$ for the sheaf of holomorphic sections of $\mathcal E$. 
Let $\Omega$ be the sheaf of differentials on $S$.  
A holomorphic connection $\nabla$ on $\mathcal E$ 
is a $\C$-linear map $\nabla: \mathcal E \to \mathcal E\otimes \Omega$ 
satisfying Leibniz rule
$\nabla(f\phi) = f\nabla(\phi) + \phi \otimes df$,  
for any holomorphic function $f$. 
\section{Appendix: The $\rho$-formalism of raising the genus of a Riemann surface}
\label{rho}
Here we recall so called $\rho$-formalism of raising the genus, 
i.e., a specific way of attaching a handle to a Riemann surface $\Sigma^{(g)}$ 
of genus $g$ to form a genus $g+1$ Riemann surface $\Sigma^{(g+1)}$ 
was introduced in \cite{Y}. 
Let $z_1$, $z_{2}$ be local  
coordinates in the neighborhood of two separated points $p_1$ and $p_2$ on $\Sigma^{(g)}$. 
Consider two disks $\left\vert z_a\right\vert \leq r_a$,
 for $r_a>0$ and $a=1,2$.
$r_1$, $r_2$ required to be small enough  so that the disks do no intersect. 
 Introduce a complex parameter $\rho$ with 
$|\rho|\leq r_1r_2$ and excise the disks
$\{z_a:\, \left\vert z_a\right\vert <|\rho |r_{\bar{a}}^{-1}\}\subset  
\Sigma^{(g)}$,  
to form a twice-punctured surface 
$\widehat{\Sigma}^{(g)}=\Sigma^{(g)}\backslash 
\bigcup_{a=1,2}\{z_a:\, \left\vert z_a\right\vert <|\rho |r_{\bar{a}}^{-1}\}$.  
We notate $\bar 1=2$, $\bar 2=1$. 
 Next define annular regions 
$\mathcal{A}_a \subset \widehat{\Sigma}^{(g)}$ with 
$\mathcal{A}_a=\{z_a:\, |\rho|r_{\bar{a}}^{-1}
\leq \left\vert z_a\right\vert \leq r_a\}$ 
and
identify them as a single region 
$\mathcal{A}=\mathcal{A}_{1}\simeq \mathcal{A}_2$ via the sewing relation 
\begin{equation}
z_1 z_2 =\rho,   
\label{rhosew}
\end{equation}
to form a compact Riemann surface 
$\Sigma^{(g+1)}=\widehat{\Sigma}^{(g)}\backslash \{\mathcal{A}_1 
\cup \mathcal{A}_2 \}\cup \mathcal{A}$
 of genus $g+1$. 
The relation (\ref{rhosew}) 
 parametrizes  a cylinder connecting the punctured Riemann surface to itself. 
On $\Sigma^{(g+1)}$ we define 
the standard basis of cycles 
$\{a_1, b_1,\ldots, a_{g+1},  b_{g+1}\}$ where the set 
$\{a_1, b_1, \ldots, a_g, b_g \}$ is the original basis on $\Sigma^{(g)}$. 
Introduce a closed anti-clockwise contour
 $\mathcal{C}_a(z_a)\subset \mathcal{A}_a$
parametrized by $z_a$ around the puncture  
at $z_a=0$.
Due to  the sewing relation \eqref{rhosew}
 $\mathcal{C}_2 (z_2)\sim -\mathcal{C}_1(z_1)$ 
 We then introduce 
 the cycle $a_{g+1}\sim \mathcal{C}_2(z_2)$ 
 and $b_{g+1}$ as 
a path chosen in $\widehat{\Sigma}^{(g)}$ 
between identified points $z_1=z_0$ and $z_2=\rho/z_0$ on the sewn surface.  
\section{Appendix: genus $g$ generalizations of elliptic functions}
\label{elliptic}
In this Appendix we recall \cite{T2} genus $g$ generalizations of classical elliptic functions. 
\subsection{Classical elliptic functions}  
\label{subsect_Elliptic}
Here we recall the classical elliptic functions \cite{Se, La}. 
For an integer $k\geq 2$, the Eisenstein series is given by  
\begin{align*}
E_k(\tau) = E_k(q) = \delta_{n, even}  
\left(-(k!)^{-1} B_k + 2((k-1)!)^{-1}\sum_{n\geq 1}\sigma_{k-1}(n)q^{n}\right),  
\end{align*}
where $\tau\in \HH$,  $q=e^{2\pi i \tau}$,   
$\sigma_{k-1}(n) = \sum_{d\vert n} d^{k-1}$, 
 and $B_k$ is the $k-{\mathrm{th}}$ 
 Bernoulli number.  
For integer $k\ge 1$, define elliptic functions $z\in \C$ 
\begin{eqnarray*}
&&
 P_1(z,\tau)=\frac{1}{z}-\sum_{k\geq 2}E_{k}(\tau )z^{k-1}, \quad 
P_k(z,\tau)=\frac{(-1)^{k-1}}{(k-1)!}\partial_z^{k-1}P_1(z,\tau),
\end{eqnarray*}
In particular 
$P_2(z, \tau) =\wp (z,\tau)+E_2(\tau)$,   
 for Weierstrass function $\wp (z,\tau)$ with periods $2\pi i$ and $2\pi i\tau$.
$P_1(z, \tau)$ is related to the quasi--periodic Weierstrass $\sigma$--function with 
$P_1(z+2\pi i\tau, \tau)=P_1(z, \tau)-1$.  
\subsection{Genus $g$ generalizations of elliptic functions}
The generalizations of elliptic functions at genus $g$ were proposed in \cite{TW}. 
Introduce a column vector 
$X=(X_a(m))$,  
 indexed by $ m\ge 0$ and $ a\in\I$ 
\[
X_a(m)=\rho_a^{-\frac{m}{2}}\sum_{\bm{b}_+}Z^{(0)}(\ldots;u(m)b_a, w_a;\ldots), 
\]
and a row vector 
$p(x)=(p_a(x,m))$,  
  for $m\ge 0, a\in\I$ 
\begin{equation*}
p_a(x,m)=\rho_a^{\frac{m}{2}}\del^{(0,m)}\psi_{p}^{(0)}(x,w_a). 
\end{equation*}
Let us also define 
column vector 
$G^{(g)}=\left(G^{(g)}_a(m)\right)$,  
 for $m\ge 0, a\in\I$, given by 
\begin{align*}
G^{(g)}=\sum_{k=1}^n \sum_{j\ge 0}\del_k^{(j)} \; q(y_{k, g})\;  
\F^{(g)}_W((u(j))_k {\bf x}_{n, g}),  
\end{align*}
where $q(y)=(q_{a}(y;m))$,  for $m\ge 0$, $a\in\I$, 
is a column vector 
\begin{equation*}
q_a (y;m)=(-1)^p \rho_a^{\frac{m+1}{2}}\del^{(m,0)}\psi_p^{(0)}(w_{-a},y),
\end{equation*} 
$R=(R_{ab}(m,n))$,  
 for $m$, $n\ge 0$ and $a$, $b\in\I$  
\begin{equation*}
R_{ab}(m,n)=\begin{cases}(-1)^p \rho_a^{\frac{m+1}{2}}
\rho_b^{\frac{n}{2}}\del^{(m,n)}\psi_{p}^{(0)}(w_{-a},w_b),& a\neq-b,
\\
(-1)^p \rho_a^{\frac{m+n+1}{2}}\E_m^n(w_{-a}),& a=-b,  
\end{cases}
\end{equation*}
\begin{equation*}
\E_m^n(y)=\sum_{\ell=0}^{2p-2}\del^{(m)}f_{\ell}(y)\;\del^{(n)}y^{\ell}, 
\quad 
\psi_p^{(0)}(x,y)=\frac{1}{x-y}+\sum_{\ell=0}^{2p-2}f_{\ell}(x)y^{\ell},
\end{equation*}
for any Laurent series $f_{\ell}(x)$ for $\ell=0,\ldots, 2p-2$. 
Define the matrices  
$\Delta_{ab}(m,n)=\delta_{m,n+2p-1}\delta_{ab}$, 
$\widetilde{R}=R\Delta$, and 
$\left(I-\widetilde{R}\right)^{-1}=\sum_{k\ge 0}\widetilde{R}^{\,k}$.  
Introduce  
$\chi(x)=(\chi_a(x;\ell))$ and  
$o(u;\bm{v,y})=(o_a(u;\bm{v,y};\ell))$, which are 
 are finite row and column vectors for  
$a\in\I$, $0\le \ell\le 2p-2$ with
\begin{equation*}
\chi_a(x;\ell)=\rho_a^{-\frac{\ell}{2}}(p(x)
+\widetilde{p}(x)(I-\widetilde{R})^{-1}R)_a(\ell), \quad 
o_a (\ell)=o_a (u;\bm{v,y};\ell)=\rho_a ^{\frac{\ell}{2}}X_a (\ell), 
\end{equation*} 
$\widetilde{p}(x)=p(x)\Delta$. 
Note that   
 $\psi_p (x,y)$ is defined by 
\begin{equation*}
\psi_p(x,y)=\psi_p^{(0)}(x,y)+\widetilde{p}(x)(I-\widetilde{R})^{-1}q(y).
\end{equation*}
For each $a \in\Ip$ introduce a vector
$\theta_a (x)=(\theta_a (x;\ell))$,    
$0\le \ell\le  2p-2$, 
\begin{equation*}
\theta_a (x;\ell) = \chi_a (x;\ell)+(-1)^p \rho_a^{p-1-\ell}\chi_{-a}(x;2p-2-\ell).
\end{equation*}
We then have the following vectors of differential forms
\begin{equation}
 P(x) =p(x) \; dx^p, \quad 
 Q(y)=q(y)\; dy^{1-p}, \quad 
\widetilde{P}(x)=P(x)\Delta, 
\end{equation}
\begin{equation}
\label{psih}
\Psi_p (x,y) =\psi_p (x,y) \;dx^p \; dy^{1-p}= 
 \Psi_p^{(0)}(x,y)+\widetilde{P}(x)(I-\widetilde{R})^{-1}Q(y). 
\end{equation}
Finally, one introduces 
\begin{equation}
\label{thetanew}
\Theta_a (x;\ell) =\theta_a (x;\ell)\; dx^p,  \quad  
O_a (u; \bm{v,y};\ell) = o_a (u; \bm{v,y};\ell) \; \bm{dy^{\wt(v)}}. 
\end{equation}
\section{Appendix: examples of vertex operator algebra $n$-point functions}
\label{examples}
\subsection{Vertex operator algebra $n$-point functions on Riemann sphere} 
\label{sphere}
For ${\bf v}_n \in V$, and a homogeneous $u \in V$, the $n$-point function on the sphere 
is given by \cite{FHL, FLM}  
\[
\F^{(0)}_W\left({\bf x}_{n, 0}\right)= 
 \langle u^{\prime }, Y(x_1) \ldots Y(x_n) u \rangle,    
\] 
while the partition function is 
$\F^{(0)}_{W, 0}=\langle u'_{(a)}, u_{(b)} \rangle= \delta_{a,b}$.   
The reductive differentials of \eqref{poros} are 
\begin{eqnarray}
\label{zhu_reduction_genus_zero_1} 
  D_1^n (x_{n+1, 0}, 0). \F^{(0)}_W\left({\bf x}_{n, 0} \right)
 &=& T_1(o(v)). \F^{(0)}_W\left({\bf x}_{n, 0} \right), 
\end{eqnarray}
\begin{eqnarray*}
  D_2^n(x_{n+1, 0}, 0). \F^{(0)}_W\left({\bf x}_{n, 0} \right) 
=z_{n+1}^{-\wt(v)} \sum\limits_{k=1}^{n} \sum\limits_{m \ge 0}    
 f_{wt(v_{n+1, 0}), m}(z_{n+1}, z_r)  T_k( v(m) ).
\F^{(0)}_W\left({\bf x}_{n, 0} \right), && 
\end{eqnarray*}
where 
 we define $f^{(0)}_{wt(v, m} (z, w)$ is a rational function defined by 
\[
f^{(0)}_{n,m}(z,w) 
=
\frac{z^{-n}}{m!}\left(\frac{d}{dw}\right)^m \frac{w^n}{z-w}, \; \;  
\iota_{z,w}f^{(0)}_{n,m}(z,w) 
=
\sum\limits_{j\in {\mathbb N}}\left( { n+j \atop m}\right) z^{-n-j-1}w^{n+j-1}, 
\]
where $\iota_{z,w}: \C[z_1, \ldots, z_n]\to \C[[z_1,z_1^{-1} \ldots, z_n z_n^{-1}]]$
 are maps \cite{FHL}.  
\subsection{Vertex operator algebra $n$-point functions on the torus}
\label{torus}
For ${\bf v}_n \in V^{\otimes n}$
 the genus one $n$-point function is defined by    
\[
\F^{(1)}_W({\bf x}_{n, 1}) 
=
Tr_W\left(Y\left(q_1^{L(0)} v_1, q_1\right) \ldots 
Y\left(q_n^{L(0)} v_n, q_n\right)  \; q^{L(0)-c/24}\right),  
\]
for $q=e^{2\pi i \tau}$ and $q_i=e^{z_i}$, where $\tau$ 
is the torus modular parametr,   
and $c$ is the central charge of the Virasoro algebra of $V$.     
For any $v_{n+1, g}\in V$, ${\bf v}_n \in V^{\otimes n}$,  
the torus reduciton formula is given by \cite{Zhu} 
\begin{eqnarray}   
 D_1^{n+1}(x_{n+1, 1}, 1).  
\F^{(1)}_W \left({\bf x}_{n, 1} \right) &=&
 \F^{(1)}_W \left( o( v_{n+1, g}) \; {\bf x}_{n,1} \right), 
\\
  D_2^{n+1}(x_{n+1, 1}, 1).
\F^{(1)}_W \left( {\bf x}_{n, 1} \right) &=&\sum\limits_{k=1}^n\sum\limits_{m \geq 0}  
 P_{m+1}
(z_{n+1}-z_k,\tau )\; 
\F^{(1)}_W ( (v[m])_k. \; {\bf x}_{n, 1}).   
\nonumber
\label{zhu_reduction_genus_one}
\end{eqnarray}
Here $P_m(z,\tau)$ denote Weierstrass functions defined by 
\[
P_m (z,\tau )=\frac{(-1)^m}{(m-1)!}\sum\limits_{n \in {\mathbb Z}_{\neq 0} }  
\frac{n^{m-1}q_z^n}{1-q^n}.
\]  
\subsection{Vertex operator algebra reduction formulas 
in genus $g$ Schottky uniformization}  
\label{corfug}
In this Section we recall reduction relations  
 for vertex operator algebra $n$-point functions  
defined on a genus $g$ Riemann surface constructed 
 in the Schottky uniformization procedure \cite{Be2, TZ, TW, T2}. 
In this case, the coefficients in reduction formulas 
are meromorphic functions on Riemann surfaces and represent genus $g$   
generalizations of the 
elliptic functions \cite{La, Se}. 
For $2g$ vertex operator algebra $V$ states 
and corresponding local coordinates 
$\bm{b}=(b_{-1}, b_1; \ldots;  b_{-g} ;b_g)$, 
$\bm{w}= (w_{-1}, w_1; \ldots; w_{-g},w_g)$,  
of $2g$ points $(p_{-1}, p_1; \ldots;  p_{-g}, p_g)$ 
on the Riemann sphere  
consider the genus zero $2g$-point correlation function
\begin{align*}
\F^{(0)}_V(\bm{b,w})=&\F^{(0)}_V(b_{-1},w_{-1};b_1, w_1;\ldots;b_{-g},w_{-g};b_g,w_g) 
\\
=&\prod_{a\in\Ip}\rho_a^{\wt(b_a)}\F^{(0)}_V 
(\bbar_1,w_{-1};b_1,w_1; \ldots; \bbar_g,w_{-g}; b_g,w_g).
\end{align*}
where $\Ip=\{1,2,\ldots,g\}$. 
Let us denote ${\bf b}_{+,g}=(b_1,\ldots,b_g)$, and 
 an element of a $V$-tensor product 
$V^{\otimes g}$-basis with the dual basis 
 ${\bf b}_{-, g}=(b_{-1}, \ldots,b_{-g})$,  
 with respect to the bilinear form 
$\langle \cdot, \cdot\rangle_{\rho_a}$ 
(cf. Appendix \ref{vertex}). 
 Let $w_a$ for $a\in\I$ be $2g$  
 formal variables and 
 ${\bm \rho}_g=(\rho_1, \ldots, \rho_g)$ $g$ complex parametrs. 
 We may identify ${\bm \rho}_g$ 
 with the canonical Schottky parametrs.   
 One introduces the genus $g$ partition function
 (zero-point function) as
\begin{align}
\label{GenusgPartition}
\F_V^{(g)} =\F_V^{(g)}(\bm{w}, {\bm \rho}_g) 
=\sum_{{\bf b}_{+,g}} \F^{(0)}_V(\bm{b, w}),   
\quad (\bm{w}, {\bm \rho}_g)=(w_{\pm 1},  \rho_1; \ldots; w_{\pm g}, \rho_g). 
\end{align}
For ${\bf x}_{n, g}=({\bf v}_{n, g}, {\bf y}_{n, g})$,  
one defines the genus $g$  formal $n$-point function for 
${\bf v}_{n, g} \in V^{\otimes n}$ and formal parametrs 
${\bf y}_{n, g}$ by 
\begin{eqnarray}
\label{GenusgnPoint}
\F_V^{(g)}({\bf x}_{n, g}) =\F_V^{(g)}({\bf x}_{n, g}; \bm{w,\rho})
=
\sum_{ {\bf b}_{+, g}} \F^{(0)}_V({\bf x}_{n, g}; \bm{b,w}), 
\nn
\F^{(0)}_V({\bf x}_{n, g}; \bm{b,w})= 
\F^{(0)}_V({\bf x}_n; b_{-1}, w_{-1};\ldots; b_g, w_g). 
\end{eqnarray}
Let $U\subset V$ be a vertex operator subalgebra such that 
 $V$ admits a $U$-module $W_\alpha$ decomposition 
$V=\bigoplus_{\alpha\in A} W_\alpha$, 
over an indexing set $A$.   
For a tensor product of $g$ modules  
$W_{\bm{\alpha}}=\bigotimes_{a=1}^g  W_{\alpha_a}$, consider 
\begin{align}
\label{eq:Z_Walpha}
\F_{W_{\bm{\alpha}}}^{(g)}({\bf x}_{n, g})  
=\sum_{{\bf b_{+, g}}\in W_{\bm{\alpha}}} \F^{(0)}_W({\bf x}_{n, g}; \bm{b,w}),  
\end{align}
where here the sum is over a basis $\{{\bf b}_{+, g}\}$ for $W_{\bm{\alpha}}$. 
It follows that
\begin{align}
\label{eq:Z_WalphaSum}
\F_W^{(g)}({\bf x}_{n, g})=\sum_{\bm{\alpha}\in\bm{A}} 
 \F_{W_{\bm{\alpha}}}^{(g)}({\bf x}_{n, g}), 
\end{align}
with  
$\bm{\alpha}=(\alpha_1,\ldots ,\alpha_g) \in \bm{A}$,  for $\bm{A}=A^{\otimes{g}}$.  
Finally, one defines corresponding formal 
$n$-point correlation differential forms 
\begin{eqnarray*} 
\widetilde{\F}_{W_{\bm{\alpha}}}^{(g)}({\bf x}_{n, g}) 
=\F_{W_{\bm{\alpha} }}^{(g)}({\bf x}_{n, g})\;
\prod_{k=1}^n dy_{k, g}^{\wt({\bf v}_{k, g})}. 
\end{eqnarray*}
Corresponding differential $D^g$ acts as 
\begin{eqnarray}
 \widetilde{\F}_V^{(g+1)}({\bf x}_{n, g})
&=& D^g. \F_V^{(g+1)}({\bf x}_{n, g}) 
= D^g. \sum_{{\bf b}_{+, g}} 
\widetilde{\F}^{(0)}({\bf x}_{n, g}; {\bf b}_{2g},{\bf w}_{2g})  
\nn
&=& \sum_{{\bf b}_{+,g+1}} \widetilde{\F}^{(0)}({\bf x}_{n, g}; {\bf b}_{2g+1},{\bf w}_{2g+1}). 
\end{eqnarray}
In \cite{TW} they prove that 
 the genus $g$ $(n+1)$-point formal differential 
$\widetilde{\F}_{W_{\bm{\alpha}}}^{(g)}(x_{n+1, g}; {\bf x}_{n, g})$, 
 for $x_{n+1, g}= ( v_{n+1, g}, y_{n+1, g})$, 
for quasiprimary vectors $v_{n+1, g} \in U$ of weight 
$\wt(v_{n+1, g})=p$ with formal parametrs 
 ${\bf y}_{n+1, g}$,     
 and general vectors ${\bf v}_n$ with parametrs 
${\bf y}_n$ satisfies the reduction formulas  
\begin{eqnarray}
\label{eq:ZhuGenusg} 
\widetilde{\F}_{ W_{ \bm{\alpha}, n+1}  }^{(g)}
\left({\bf x}_{n+1, g}  
\right) &=& \left(D_1^{(n+1)} + D_2^{(n+2)}\right). 
\widetilde{\F}^{(g)}_{W_{\bm{\alpha}}, n } \left({\bf x}_{n, g} \right), 
\end{eqnarray}
\begin{eqnarray*}
D_1^{n+1} ( x_{n, g}, g). \widetilde{\F}^{(g)}_{W_{\bm{\alpha}}, n }
\left( {\bf x}_{n, g} \right) 
= \sum_{a=1}^g \Theta_a (y_{n+1, g}) \; O_a^{ W_{ \bm{\alpha} }}  \;  
\left(v_{n+1, g}; {\bf x}_{n, g}\right), \qquad \qquad \qquad \qquad &&
\nn
D_2^{n+1}(x_{n+1, g}, g).\widetilde{\F}^{(g)}_{ W_{\bm{\alpha}}, n+1}
\left({\bf x}_{n+1, g}  \right)= \sum_{k=1}^n \sum_{j\ge 0}\del^{(0,j)}  
  \Psi_p(y_{n+1, g},y_{k, g})\;  
 \widetilde{\F}_{W_{\bm{\alpha}}, n}^{(g)}  
\left( (u(j))_k. {\bf x}_{n, g} \right) dy_{k, g}^j, &&  
\end{eqnarray*}
Here $\del^{(0,j)}$ is given by     
$\del^{(i,j)}f(x,y)=\del_x^{(i)}\del_y^{(j)}f(x,y)$,   
for a function $f(x,y)$, and $\del^{(0,j)}$ 
denotes partial derivatives with respect to $x$ and $y_{j, g}$.  
The forms 
$\Psi_p(y_{n+1, g}, y_{k, g} )\; dy_{k, g}^j$ are given by   
 \eqref{psih}, 
 $\Theta_a (x)$ is of \eqref{thetanew}, and 
 $O^{ W_{\bm{\alpha}}}_a (v_{n+1, g}; {\bf x}_{n, g})$ 
is \eqref{thetanew}.   

\end{document}